% Gromov-Witten invariants of flag manifolds, via D-modules
% A. Amarzaya and M. A. Guest
% This is an AMS-TeX file

\input amstex
\documentstyle{amsppt}

\magnification=1200
\pagewidth{5.4in}
\pageheight{7.5in}
\parskip 10pt

\expandafter\redefine\csname logo\string@\endcsname{}
\NoBlackBoxes                
\NoRunningHeads
\redefine\no{\noindent}

\define\C{\bold C}

\define\Z{\bold Z}

\define\al{\alpha}
\define\be{\beta} 

\define\de{\delta}
\define\la{\lambda}

\define\La{\Lambda}

\define\Om{\Omega}

\define\th{\theta}
\define\om{\omega}

\define\thi{\theta^{(i)}}
\define\thj{\theta^{(j)}}

\define\thzero{\theta^{(0)}}
\define\thone{\theta^{(1)}}
\define\thtwo{\theta^{(2)}}
\define\ththree{\theta^{(3)}}
\define\thfour{\theta^{(4)}}

\define\st{\ \vert\ }   
\redefine\ll{\lq\lq}
\redefine\rr{\rq\rq\ }
\define\rrr{\rq\rq}

\redefine\det{\operatorname {det}}
\redefine\deg{\operatorname {deg}}
\redefine\dim{\operatorname {dim}}

\define\rank{\operatorname {rank}}

\define\glnc{GL_n\C}
\define\glsc{GL_{s+1}\C}

\define\Gl{GL}

\redefine\b{\partial}

\topmatter

\title Gromov-Witten invariants of flag manifolds, via D-modules
\endtitle

\author A. Amarzaya and M. A. Guest
\endauthor

\thanks The authors were partially supported by a 
grant from the JSPS.
\endthanks

\endtopmatter

\document

We  present a method for computing the $3$-point genus zero Gromov-Witten
invariants of the complex flag manifold $G/B$
from the relations of the small quantum cohomology algebra $QH^\ast G/B$ ($G$ is a
complex semisimple Lie group and $B$ is a Borel subgroup). 
In \cite{Fo-Ge-Po} and \cite{Ki-Ma}, at least in
the case $G=\glnc$, two algebraic/combinatoric methods have been proposed,
based on suitably designed axioms.
Our method is quite different,
being differential geometric in nature; it is based on the approach to quantum
cohomology described in \cite{Gu}, which is in turn based on the integrable
systems point of view of Dubrovin and Givental.

In \S 1 we shall review briefly the method of \cite{Gu}.   In \S 2 we discuss the
special properties of $G/B$ which lead to a computational algorithm. In fact the
same method works for any Fano manifold whose cohomology is generated
by two-dimensional classes, so our approach is
more general than those of \cite{Fo-Ge-Po} and \cite{Ki-Ma}.
In \S 3 we present explicit results for the
case $G=\glnc$, $n=2,3,4$.  In \S 4, we show how to produce
\ll quantum Schubert polynomials\rr for $G/B$, by which we mean specific polynomial representatives of quantum Schubert classes.

The second author thanks Josef Dorfmeister for essential suggestions concerning the proof of
Proposition 2.2. We are very grateful to the referee for pointing out several
inaccuracies in an earlier version of this paper.

\head
\S 1 Quantum cohomology via D-modules
\endhead

We list some well known properties of the cohomology and
quantum cohomology algebras of $G/B$ (see 
\cite{Ci}, \cite{Gi-Ki}, \cite{Ki}).
The cohomology algebra of $G/B$ (with complex coefficients) 
has the form
$$
H^\ast G/B = \C[b_1,\dots,b_r]/(R_1,\dots,R_u),
$$
where $b_1,\dots,b_r$ are additive generators of $H^2 M$ (with $r=\rank G$), and
$R_1,\dots,R_u$ are explicitly known homogeneous polynomial relations in $b_1,\dots,b_r$.
The (small) quantum cohomology
algebra of $G/B$ is of the form
$$
QH^\ast G/B = \C[b_1,\dots,b_r,q_1,\dots,q_r]/(\Cal R_1,\dots,\Cal R_u),
$$
where
$\Cal R_1,\dots,\Cal R_u$ are explicitly known polynomial relations in $b_1,\dots,b_r$
and additional variables $q_1,\dots,q_r$, with the property  $\Cal R_i\vert_{q=0}=R_i$
for all $i$. The relations $\Cal R_i$ are homogeneous (in the sense
that all terms in the same relation have the same degree) if $b_1,\dots,b_r$ have their
usual degrees and  $q_1,\dots,q_r$ have  degree $4$.  Moreover,
$H^\ast G/B\otimes \C[q_1,\dots,q_r]$ and $QH^\ast G/B$ are additively isomorphic  (they
are both free modules of rank $s+1=\dim H^\ast G/B$ 
over the polynomial ring $\C[q_1,\dots,q_r]$)
but they are not isomorphic as 
$\C[q_1,\dots,q_r]$-algebras.
This is the starting point for \cite{Fo-Ge-Po} and \cite{Ki-Ma} as well as 
for us.  

The objective is to compute the $3$-point genus zero Gromov-Witten invariants of $G/B$, or, equivalently,  the quantum product operation on  $H^\ast G/B\otimes \C[q_1,\dots,q_r]$.  This, in turn, is equivalent to a particular choice of additive isomorphism
$$
QH^\ast G/B\to H^\ast G/B\otimes \C[q_1,\dots,q_r].
$$
Our strategy is to characterize this isomorphism,  and then to compute it explicitly.

To begin, we define an (abstract) {\it quantum evaluation map} as a map
$\de:QH^\ast G/B\to H^\ast G/B\otimes \C[q_1,\dots,q_r]$
with the following four properties.  In stating these properties,
for a polynomial $c$ in the variables $b_1,\dots,b_r,q_1,\dots,q_r$
we denote the corresponding element of $QH^\ast G/B$ by $[c]$,
and the corresponding element of 
$H^\ast G/B\otimes \C[q_1,\dots,q_r]$ by $[[c]]$. 

\no($\de$1) $\de$ is an isomorphism of $\C[q_1,\dots,q_r]$-modules,

\no($\de$2) $\de[1]=[[1]]$ and  $\de[b_i]=[[b_i]]$ for $i=1,\dots,r$,

\no($\de$3)  $\deg \de[c] = \deg [c]$ for any homogeneous polynomial $c$, and

\no($\de$4)  $(\de[c])\vert_{q=0} = [[c\vert_{q=0}]]$
for any polynomial $c$.

\no{\it Example 1:}  Denote the quantum product on
$H^\ast G/B\otimes \C[q_1,\dots,q_r]$ by $\circ$.   Then we
have a quantum evaluation map $\de$ which \ll evaluates
quantum  polynomials\rrr, i.e.\   $\de[c]=c^\circ$, 
where $c^\circ$ is the element of $H^\ast G/B\otimes \C[q_1,\dots,q_r]$
obtained by replacing each monomial $b_ib_jb_k\dots$
in $c$ by $[[b_i]]\circ[[b_j]]\circ[[b_k]]\circ\dots$  (this example
motivates our terminology \ll quantum evaluation map\rrr).

\no{\it Example 2:}  
Choose a basis $[c_0],\dots,[c_s]$
of the $\C[q_1,\dots,q_r]$-module $QH^\ast G/B$, with the property that
$[[c_0\vert_{q=0}]],\dots,[[c_s\vert_{q=0}]]$ is a basis of
the $\C[q_1,\dots,q_r]$-module $H^\ast G/B\otimes\C[q_1,\dots,q_r]$.  
Assume further that the $c_i$ are homogeneous polynomials such that
$c_0=1$ and $c_i=b_i$ for $i=1,\dots,r$.  Then 
$$
\de[c_i] = [[c_i\vert_{q=0}]]
$$
defines a quantum evaluation map.

An (abstract) quantum evaluation map gives
rise to an (abstract) quantum product operation $\circ$ by the formula 
$$
x\circ y  = \de(\de^{-1}( x)  \de^{-1} (y )).
$$  
In example 1, this is the usual quantum product. In example 2, however, we
may or may not obtain the usual quantum product.  

Let us examine the
properties of the product $\circ$ from example 2 in more detail.
First, we claim that $\circ$ satisfies the following three properties (all
of which are well known properties of the usual quantum product):

\no(1) $(H^\ast G/B\otimes \C[q_1,\dots,q_r], \ \circ\ )$ is an algebra, isomorphic to the
algebra $QH^\ast G/B$.

\no(2) for any $x,y\in H^\ast G/B$, one has
$x\circ y=xy + \text{terms involving $q_1,\dots,q_r$}$, 
where $x y$ denotes the cup product of $x$ and $y$, and both
sides of this formula are homogeneous of degree $\deg x + \deg y$.

\no(3) $\de[c]=c^\circ$.

\no Property (1) is true by the definition of $\circ$.  
For (2) it suffices to check that 
$$
[[c_i\vert_{q=0}]]\circ[[c_j\vert_{q=0}]]\vert_{q=0} = 
[[c_i\vert_{q=0}]]\,[[c_j\vert_{q=0}]].
$$ 
We have
$[[c_i\vert_{q=0}]]\circ[[c_j\vert_{q=0}]]\ \vert_{q=0} = 
\de([c_i][c_j])\vert_{q=0}$ by definition,
and this is $[[c_i][[c_j]]\vert_{q=0}$ by $(\de4)$, as required.
As for (3), it suffices to observe that
$$
\align
(b_{i}b_{j}b_{k}\dots)^\circ&=
[[b_{i}]]\circ [[b_{j}]]\circ [[b_{k}]]\circ\dots \\
&= 
\de( \de^{-1}[[b_{i}]] \de^{-1}[[b_{j}]]  \de^{-1}[[b_{k}]]\dots) \ \text{by definition}
\\
&= \de( [b_{i}] [b_{j}] [b_{k}]\dots) \ \text{by property $(\de2)$}
\\
&= \de ([b_{i}b_{j}b_{k}\dots]).
\endalign
$$

The usual quantum product satisfies a further condition, however, and this is an
\ll integrability condition\rrr, which is more subtle than the
previous algebraic conditions, and it leads to the appearance of differential 
equations.  It is at this point that we diverge from \cite{Fo-Ge-Po} and \cite{Ki-Ma}
in an essential way, as those authors impose further algebraic conditions in
order to characterize the usual quantum product.

Let us introduce matrix functions $\om_1,\dots,\om_r$ of $q_1,\dots,q_r$  as follows:   
$$
[b_ic_j]=\sum_{k=0}^s (\om_i)_{kj}[c_k].
$$
Thus, $\om_i$ can be interpreted as the matrix of multiplication by $[b_i]$ on $QH^\ast G/B$ with respect to the basis $[c_0],\dots,[c_s]$, or the 
matrix of the operator $[[b_i]]\circ$ on $H^\ast G/B\otimes\C[q_1,\dots,q_r]$ 
with respect to the basis $\de[c_0],\dots,\de[c_s]$.  The integrability condition is
$$
q_i \tfrac{\b}{\b q_i} \om_j = 
q_j\tfrac{\b}{\b q_j} \om_i
\ \text{for all}\ 1\le i,j\le r.\tag$\ast$
$$
This may be written $d\om=0$, where 
$\om = \sum_{i=1}^r \om_i dt_i$ and  $q_i=e^{t_i}$.   
Since the commutativity and associativity of $\circ$ imply that
$\om\wedge\om = 0$, the integrability condition is equivalent to
$d\om +\om\wedge\om=0$.  It is usual
in quantum cohomology theory to regard 
$t_1,\dots,t_r$ as coordinates on the vector
space $H^2 G/B\cong \C^r$; then 
$d+\om$ is a connection in the trivial vector bundle
with fibre $H^\ast G/B \cong \C^{s+1}$ over the manifold
$H^2 G/B \cong \C^{r}$, and $d\om +\om\wedge\om=0$ means
that this connection has zero curvature.

It is a nontrivial matter to construct a quantum evaluation map $\de$
satisfying condition $(\ast)$.  For the usual quantum product (example 1 above), $(\ast)$ is
a consequence of the moduli space construction of quantum cohomology.  
For a product obtained from example 2, $(\ast)$ imposes a strong condition
on the choice of basis. To show that a suitable basis exists (and
that it produces the usual quantum product), 
we shall use the
method of \cite{Gu}, which depends on the existence of a
{\it quantization} of the algebra $QH^\ast G/B$. To define this concept,
we work in the algebra $D^h$ of  differential operators 
generated by $h\b_1,\dots,h\b_r$ with coefficients in 
$\C[q_1,\dots,q_r,h]$, 
where $h$
is a complex parameter, and where $\b_i = \frac{\b}{\b t_i} = q_i\frac{\b}{\b q_i}$.

\proclaim{Definition 1.1} A {\it quantization} of the algebra $QH^\ast G/B$ is a (left, cyclic)
$D^h$-module $M^h= D^h/(D_1^h,\dots,D_u^h)$ such that

\no (1) $M^h$ is free over $\C[q_1,\dots,q_r,h]$
 of rank $s+1$ (where $s+1=\dim H^\ast G/B$),

\no (2) $\lim_{h\to 0} S(D^h_i) = \Cal R_i$,
where $S(D^h_i)$ is the result of replacing $h\b_1,\dots,h\b_r$
by $b_1,\dots,b_r$ in $D^h_i$ (for $i=1,\dots,u$).
\endproclaim

A quantization  of $QH^\ast G/B$ is provided by the so-called
quantum Toda Lattice associated to the Langlands dual of $G$ (this can be proved
without any reference to quantum cohomology theory --- see \cite{Go-Wa}
and \cite{Ko}). 
Let $P_0,\dots,P_s$ to be the \ll standard monomials\rr 
in $h\b_1,\dots,h\b_r$ with respect to a choice
of Gr\"obner basis for the ideal $(D_1^h,\dots,D_u^h)$;
the equivalence classes  $[P_0],\dots,[P_s]$
form an additive basis for $M^h$.  Let
$c_i=S(P_i)$ where the notation $S$ is 
as in (2) above; then each $c_i$ is a monomial in $b_1,\dots,b_r$
and $[c_0],\dots,[c_s]$ is a basis of $QH^\ast G/B$. 

From $c_0,\dots,c_s$ we obtain
a quantum evaluation map $\de$ satisfying $(\de1)-(\de4)$, in
the manner explained earlier.
However, in general, this $\de$ does not satisfy the integrability condition $(\ast)$.
We need modified polynomials $\hat c_0,\dots,\hat c_s$ (with associated
quantum evaluation map
$\hat\de$), and it is this modification which makes essential
use of the underlying $D^h$-module $M^h$.

Let $\Om^h_i$
denote the matrix of the action
of $\b_i$ on $M^h$ with respect to the basis
$[P_0],\dots,[P_s]$; more precisely we define the functions
$h\Om^h_1,\dots,h\Om^h_r$ by
$$
[h\b_i P_j]=\sum_{k=0}^s (h\Om^h_i)_{kj}[P_k].
$$
Let $\Om^h=\sum_{i=1}^r \Om^h_i dt_i$. 
Then $\Om^h$ is of the form
$$
\Om^h = \frac 1h \om + \th^{(0)} + h\th^{(1)} + \dots + h^p \th^{(p)}
$$
where $\th^{(0)},\dots,\th^{(p)}$ depend polynomially on $q_1,\dots,q_r$.

As in the case of $\om$, it is convenient to regard $d+\Om^h$
as a (family of) connections in the trivial vector bundle
with fibre $H^\ast G/B \cong \C^{s+1}$ over the manifold
$H^2 G/B \cong \C^{r}$.  This connection is flat for every $h$ 
(see \cite{Gu}, \S 1).
If all $\th^{(i)}$ were
zero we would be able to conclude that $d+\om$ is flat, and that 
condition $(\ast)$ is satisfied.
We shall obtain a modified $\hat\Om^h$ which has this property.

Let us restrict the parameter $h$ to the unit circle 
$S^1 = \{ h\in\C \st \vert h\vert = 1 \}$.
Since $d + \Om^h$ is flat, and $\C^r$ is simply connected,
we have $\Om^h = L^{-1}dL$
for some $L:\C^r\to \La \glsc$, where $\La \glsc$ is the
(smooth) loop group of $\glsc$, i.e.\  the space of all (smooth)
maps $S^1\to \glsc$.  
At this point we leave (temporarily) the polynomial algebra
$\C[q_1,\dots,q_r,h]$;  $L$ is not in general a polynomial
function of $q_1,\dots,q_r,h$.  
Let $L=L_-L_+$ be the Birkhoff factorization
(see \cite{Pr-Se}, Chapter 8) of $L$ at $q=q_0$, with
$$
\align
L_-(q,h)&=I+h^{-1}A_1(q)+h^{-2}A_2(q)+\dots\\
L_+(q,h)&=B_0(q)+hB_1(q)+h^2B_2(q)+\dots
\endalign
$$
(Without loss of generality we can assume that $L\vert_{q_0}$ admits such
a factorization for some point $q_0$, hence the same is true of $L\vert_{q}$ for $q$ in
a neighbourhood of $q_0$.)
The gauge transformation
$L\mapsto \hat L = L(L_+)^{-1}=L_-$ transforms $\Om^h = L^{-1}dL$
into $\hat\Om^h =  L_-^{-1} dL_{-}$, and a simple
calculation (see \cite{Gu}, \S 3) shows that 
$$
\hat \Om^h =\frac 1h\hat\om
$$ 
where $\hat\om=  B_0 \om B_0^{-1}$.  
We shall show, in the next section, that $L_+$ is polynomial in $q_1,\dots,q_r,h$.
In particular it is regular at $q=q_0=0$, and we may normalize it by 
insisting that $L\vert_{q=0}=I$, which determines $L_+$ uniquely.

Note that $\hat \Om^h_i$ is
the matrix of the action of $\b_i$ with respect to
$[\hat P_0],\dots,[\hat P_s]$, where
$\hat P_i = \sum_{j=0}^s (L_+)^{-1}_{ji} P_j$,
and $\hat\om_i$ is
the matrix of multiplication by $[b_i]$
with respect to
$[\hat c_0],\dots,[\hat c_s]$, where
$\hat c_i = \sum_{j=0}^s (B_0^{-1})_{ji} c_j$.

The connection $d+\hat \Om^h$, hence also $d+\hat\om$, is
 flat, so we have:

\proclaim{Theorem 1.2} The quantum evaluation map $\hat\de$ defined using the
polynomials $\hat c_0,\dots,\hat c_s$ satisfies the integrability condition $(\ast)$.
\qed\endproclaim

The same 
$D^h$-module provides a quantization,
in the above sense,  when $\circ$
is the usual quantum product (see \cite{Ki}). 
We shall use this, in the next section,
to prove:

\proclaim{Theorem 1.3} Our product operation
(obtained from the quantum evaluation map $\hat\de$) agrees with
the  usual quantum product operation. 
\endproclaim

In the rest of the article we shall set up an algorithm for computing
the quantum product. The procedure  may be summarized
as follows:  first we choose $P_0,\dots,P_s$, then we compute $\Om^h$
and $L_+$, and finally  the constant term of $L_+$ produces 
$\hat c_0,\dots,\hat c_s$.  The quantum evaluation map is then
given explicitly by
$$
\hat c_i^\circ = \hat\de[\hat c_i] = [[\hat c_i\vert_{q=0}]]=[[c_i]].
$$
These formulae, together with the relations of the quantum cohomology algebra,
determine the quantum product completely.

Another way of computing the quantum product would be to read off the products of degree two classes with arbitrary classes from the matrices $\hat\om_i$.  From the above description of $\hat\om_i$, it is 
the matrix of the operator $[[b_i]]\circ$
with respect to
$\hat\de[\hat c_0],\dots,\hat\de[\hat c_s]$. Since
$\hat\de[\hat c_j]=[[\hat c_j\vert_{q=0}]]=[[c_j]]$, this
basis is just the original basis $[[c_0]],\dots,[[c_s]]$.  Thus
we obtain all products $[[b_i]]\circ [[c_j]]$ for
$1\le i\le r$, $0\le j\le s$, and hence all other products.

The main difficulty in the algorithm is the computation of $L_+$, and the next section will be devoted to this.

\head
\S 2 The p.d.e.\  system for $L_+$
\endhead

Let us choose the bases $[P_0],\dots,[P_s]$ of $M^h$ and $[c_0],\dots,[c_s]$ 
of $QH^\ast G/B$, as explained earlier.  We obtain 
$\Om^h = \frac 1h \om + \th^{(0)} + h\th^{(1)} + \dots + h^p \th^{(p)}$,
hence there exist $L$, 
$L_-$, $L_+$, such that
$$
\frac1h B_0\om B_0^{-1} = (L_-)^{-1}dL_-  = L_+\Om^h (L_+)^{-1}  + L_+d((L_+)^{-1}).
$$
Our task, therefore,  is to solve the system of partial differential equations
$$
\frac1h B_0\om B_0^{-1} L_+ \ =\ L_+ \Om^h \ -\ dL_+\tag$\ast\ast$
$$
for $L_+(q,h)=
B_0(q)+hB_1(q)+h^2B_2(q)+\dots$ in terms of the
known $1$-forms $\Om^h$, $\om$.  This should be regarded as
an explicit form of the integrability condition $(\ast)$ of \S 1.  
Note that  only the first unknown
function $B_0$ is needed for the computation of the quantum product.  However,
{\it a priori}, the system involves infinitely many unknown functions
$B_0,B_1,\dots$ and is rather complicated.  {\it The main observation of this paper
is that the system can be solved \ll by quadrature\rrr, and that an algorithm for
the solution can be implemented by computer.}

For the computation we shall rewrite $L_+$ as
$$
L_+(q,h)=Q_0(q)(I+hQ_1(q)+h^2Q_2(q)+\dots).
$$
Comparing coefficients of powers of $h$ in $(\ast\ast)$ gives the system
$$
\align
Q_0^{-1}dQ_0 &= \thzero + [Q_1,\om]\\
dQ_1&=\thone 
+ [Q_1,\thzero] + [Q_{2},\om]  - [Q_1,\om]Q_1\\
dQ_i&=\thi + Q_1 \theta^{(i-1)} + \dots + Q_{i-1} \thone
+ [Q_i,\thzero] + [Q_{i+1},\om]  - [Q_1,\om]Q_i\  (i\ge 2)
\endalign
$$

To proceed further we need some more notation.

\proclaim{Definition 2.1} Let $A=(a_{ij})$ be an $(s+1)\times (s+1)$ matrix.  Let $k$ be
an integer with $-s\le k\le s$.  

\no(1) The
$k$-diagonal of $A$ is the matrix whose $(i,j)$-th entry is $a_{ij}$
when $j=i+k$ and zero when $j\ne i+k$. 

\no(2) We say that $A$ is $l$-triangular if all the $k$-diagonals
of $A$ with $-s\le k < l$ are zero.  

\no (Thus the $0$-diagonal is the
usual diagonal; $A$ is $0$-triangular if
and only if it is upper triangular; $A$ is $1$-triangular if
and only if it is strictly upper triangular.)
\endproclaim

\no The matrices in $(\ast\ast)$ are $(s+1)\times (s+1)$ complex matrices, where
$s+1=\dim_{\C} H^\ast G/B$ (cohomology with complex coefficients),
and they have a natural \ll block\rr structure
$A=(A_{\al,\be})$ where  $A_{\al,\be}$ is an $s_{\al}\times s_{\be}$ matrix 
(a  block) with  $s_{\al}=\dim_{\C}  H^{2\al} G/B$,  $0\le \al \le m = \dim_{\C} G/B$.  We
shall generally use Greek indices, separated by commas, in reference
to block matrices.
We may extend Definition 2.1 in an obvious way:
we say $A$ is a (block) $l$-triangular matrix if
$A_{\al,\be}$ is a zero matrix whenever $\be<\al+l$.  {\it As we always deal with block
matrices, and always of the same form, we shall omit the word \ll block\rr from
now on.}

By Proposition 3.3 of \cite{Gu} we have the following
homogeneity conditions:

\no(H1) each entry of $(\om_i)_{\al,\be}$ is (either identically zero or)
homogeneous of degree $2(\be-\al+1)$,

\no(H2) each entry of $(\theta_i^{(j)})_{\al,\be} $ is (either identically zero or)
homogeneous of degree $2(\be-\al-j)$.

\no Since $\Om^h$ is polynomial in $q_1,\dots,q_r$
and each $q_i$ has degree $4$,   (H1) and (H2) imply:

\no(F1) $\om$ is $(-1)$-triangular,

\no(F2) $\thj$ is $(j+2)$-triangular for $j=0,1,\dots,m-2$.

\no We shall use these properties to prove:

\proclaim{Proposition 2.2} We have

\no(1) $Q_0=\exp X$ where $X$ is $2$-triangular,

\no(2) for $i\ge 1$, $Q_i$ is $(i+2)$-triangular, and

\no(3) for $i\ge 0$, each entry of $(Q_i)_{\al,\be}$ is (either identically zero or) homogeneous
of degree $2(\be-\al-i)$.

\no In particular, $L_+=Q_0(I + hQ_1 + \dots + h^{m-2}Q_{m-2})$, where
$m=\dim_{\C} G/B$, i.e.\  there are
only finitely many nonzero matrices $Q_i$. (Recall that we have
$(m+1)\times (m+1)$ block matrices, and that \ll triangular\rr means
\ll block triangular\rrr.)
\endproclaim

\demo{Proof} Properties (F1) and (F2) show that $\Om^h$ is a $1$-form
taking values in the subspace of the Lie algebra of the loop group of $\Gl_{s+1}\C$
consisting of elements of the form
$\sum_{i\in\Z} A_i h^i$ such that $A_{i}$ is $i$-triangular for $i < 0$ and 
$(i+2)$-triangular for  $i\ge 0$ (in particular
$A_i=0$ for $i\ge m-1$).
This subspace is a Lie subalgebra
(since the product of an $a$-triangular matrix with a $b$-triangular
matrix is $(a+b)$-triangular).  Hence $L$ (where $L^{-1}dL=\Om^h$) takes
values in the corresponding subgroup of the loop group of $\Gl_{s+1}\C$,
and so do $L_-$ and $L_+$.   We have not specified the topology of the loop group
here because all arguments may be carried out in a finite dimensional quotient
group, but for the sake of definiteness one may take the smooth loop group,
where the Birkhoff decomposition is known from \cite{Pr-Se}. In particular
$(L_+)^{-1}dL_+$ takes values in the finite dimensional subalgebra consisting of 
elements of the form $\sum_{i=0}^{m-2} A_i h^i$ such that $A_i$ is $(i+2)$-triangular,
and $L_+$ takes values in the finite dimensional subgroup consisting of 
elements of the form $\sum_{i=0}^{m-2} B_i h^i$ such that $B_0=\exp X$ 
where $X$ is $2$-triangular, and
$B_i$ is $(i+2)$-triangular for $i=1,\dots,m-2$.  (This is a subgroup of a
slightly larger finite dimensional group which was described in detail in \cite{Bu-Gu}.)
Properties (1) and (2) follow immediately, and property (3) is a consequence of the
homogeneity properties (H1) and (H2) as in the proof of Proposition 3.3 of \cite{Gu}.
\qed\enddemo

To proceed further, we write
$$
\align
Q_0 &= I+Q_0^{[2]}+Q_0^{[3]}+\dots+Q_0^{[m]}\\
Q_i &= Q_i^{[i+2]}+Q_i^{[i+3]}+\dots+Q_i^{[m]}\quad(1\le i\le m-2)\\
\om &= \om^{[-1]}+\om^{[0]}+\om^{[1]}+\dots+\om^{[m]}\\
\thi &= \theta^{(i),[i+2]}+\theta^{(i),[i+3]}+\dots+\theta^{(i),[m]}\quad(1\le i\le m-2)
\endalign
$$
for the respective decompositions into diagonal parts, and then
decompose the system of equations accordingly, using
$$
(XY)^{[j]} = \sum_{k} X^{[j-k]}Y^{[k]}, \quad
[X,Y]^{[j]} = \sum_{k} [X^{[j-k]},Y^{[k]}]
$$
where the suffix $[j]$ always refers to the $j$-diagonal part.
The equation for $dQ_i$ becomes:
$$
\align
dQ_i^{[j]}=\theta^{(i),[j]} 
+
&\tsize\sum_{k=i+1}^{j-3} Q_1^{[j-k]} \theta^{(i-1),[k]}
\\
+
&\quad \dots 
\\
+
&\tsize\sum_{k=3}^{j-i-1} Q_{i-1}^{[j-k]} \theta^{(1),[k]}
\\
+
&\tsize\sum_{k=2}^{j-i-2} [Q_i^{[j-k]}, \theta^{(0),[k]}]
\\
+ 
&\tsize\sum_{k=-1}^{j-i-3} [Q_{i+1}^{[j-k]}, \om^{[k]}]
\\
-
&\tsize\sum_{k=-1}^{m-1} \tsize\sum_{l=3}^{m}
 [Q_{1}^{[l]} ,\om^{[k]}]Q_{i}^{[j-k-l]}
\endalign
$$
(for $2\le i\le m-2$),
and the equations for $Q_0$ and $Q_1$ decompose in a similar way.

\proclaim{Definition 2.3} We define a total ordering on the (symbols) 
$Q_{i}^{[j]}$,  $2\le i+2\le j\le m$ by:
$Q_{i_1}^{[j_1]} < Q_{i_2}^{[j_2]}$ if and only if 
(a) $j_1-i_1 <
j_2-i_2$ or 
(b)$j_1-i_1 =
j_2-i_2$ and $j_2<j_1$.
\endproclaim

The key observation concerning the above system is that it is of the form
$$
dQ_i^{[j]}=\ \text{an expression involving}\ Q_{l}^{[k]}
\ \text{with}\ Q_{l}^{[k]} <  Q_i^{[j]} 
$$
Starting with the smallest term $Q_{m-2}^{[m]}$, we may therefore
integrate successively to find all $Q_i^{[j]}$.  The integrability condition
(for solving $dX=Y$, i.e.\  $dY=0$) is satisfied at each step because the
Birkhoff decomposition guarantees the consistency of the system.

Finally we shall give the proof of Theorem 1.3 and the proof of the polynomiality of $L_+$:

\demo{Proof of Theorem 1.3}
 Let $P_i, c_i$ and $\hat P_i, \hat c_i$ be as in
\S 1.  Let $T_i:M^h\to M^h$ be the  map given by the action of $\b_i$. 
The matrix of $T_i$ (in the sense of \S 1) is  $\Om^h_i$, with respect to the
basis $[P_0],\dots,[P_s]$.
The procedure of \S 1 gives a quantum evaluation map
$\hat\de$ and a product operation $\circ$ such that
$\hat c_i^\circ = \hat\de [\hat c_i] = [[\hat c_i\vert_{q=0}]] = [[c_i]]$,
and we wish to show that this is the usual quantum product. 

Let $\ast$ denote the usual quantum product.  It is known that for each $i$ there exists
a (homogeneous) polynomial $\bar c_i$ such that 
$\bar c_i^\ast = [[c_i]]$.  
The usual quantum cohomology D-module $\bar M^h$ 
(see \cite{Gi}, \cite{Ki})
produces a (homogeneous) differential operator $\bar P_i$,
with polynomial coefficients, such that
$\lim_{h\to 0} S(\bar P_i) = \bar c_i$ (where $S$ is as in Definition 1.1).  These give
a basis of $\bar M^h$.  Let $\bar T_i:\bar M^h\to \bar M^h$ be the  map
given by the action of $\b_i$, and let $\bar \Om^h_i$ be 
the matrix of $\bar T_i$  with respect to the
basis $[\bar P_0],\dots,[\bar P_s]$.  It is known that $\bar\Om^h_i = \frac 1 h \bar\om$
for some $\bar\om$.

Now we use the fact (\cite{Ki}) that $M^h=\bar M^h$, hence $T_i=\bar T_i$.
Since  $P_0,\dots,P_s$ are monomials in $h\b_1,\dots,h\b_r$, we can write
$\bar P_i = \sum_{j=0}^s U_{ji} P_j$ for a (homogeneous) map
$U(q,h) 
= U_0(q) + h U_1(q) + h^2 U_2(q) + \dots$ (which is polynomial in
$q_1,\dots,q_r,h$).  The matrix of $\bar T_i$
with respect to the
basis $[P_0],\dots,[P_s]$ is therefore $U\bar\Om^h U^{-1} + UdU^{-1}$, so we
have
$\Om^h = U\bar\Om^h U^{-1} + UdU^{-1}$.
The modified connection form 
$\hat \Om^h_i$ (from \S 1) is related to $\Om^h_i$ by
$\hat\Om^h =  (L_+)^{-1}\Om^h L_+  +  (L_+)^{-1} d L_+$.
It follows that
$\hat \Om^h$ ($=\frac 1 h \hat\om$) and $\bar \Om^h$ ($=\frac 1 h \bar\om$)
differ by the (homogeneous) gauge transformation $L_+ U$.

The method used in this section to solve the p.d.e.\  $(\ast\ast)$ for $L_+$ applies also
to the analogous p.d.e.\  for $L_+U$ (with all $\th^{(i)}=0$)
and shows that $L_+ U = I$.  It follows that $L_+$ and $(L_+)^{-1}$ are polynomial
in $q_1,\dots,q_r,h$, and that
$$
\bar c_i = \lim_{h\to 0} S(\bar P_i) =
\lim_{h\to 0} S(\sum_{j=0}^s U_{ji} P_j) = 
\lim_{h\to 0} S(\sum_{j=0}^s (L_+^{-1})_{ji} P_j) = 
\sum_{j=0}^s (Q_0^{-1})_{ji} c_j = \hat c_i.
$$
This shows that $\ast$ coincides with $\circ$.
\qed\enddemo

\head
\S 3 Example: $G=\glnc$
\endhead

For $G=\glnc$ we have $r=n-1$, $m=\frac12n(n-1)$.  The dimensions  $s_0,\dots,s_m$
are determined by the Poincar\acuteaccent e polynomial
$$
\sum_{i=0}^m  s_i z^{2i} =
(1+z^2)(1+z^2+z^4)\dots(1+z^2+\dots+z^{2n-2}),
$$
and we have $s+1=n!=s_0+\dots+s_m$.

Many of the terms in our p.d.e.\  system vanish:

\proclaim{Proposition 3.1}
(1) $\om^{[j]} =0\ \ \text{if $j$ is even}$;
(2) $\theta^{(i),[j]}=0\ \ \text{if $j-i$ is odd}$.
\endproclaim

\demo{Proof} Property (H1) of \S 2 says that each entry of $\om^{[j]}$ 
is either identically zero or homogeneous of degree $2j+2$. As $\om^{[j]}$ 
is a polynomial function of $q_1,\dots,q_{n-1}$, and
$\deg q_1=\dots=\deg q_{n-1}=4$,
$\om^{[j]}$  must be  identically zero if $j$ is even.  
A similar argument
applies to  $\theta^{(i),[j]}$, using (H2).
\qed\enddemo

We will use this fact without comment from now on.  
Our strategy for solving the system and finding the quantum product has
three steps.  

\no{\it Step 1:}  Write down the system, removing terms which are known {\it a priori} to
vanish.

\no{\it Step 2:}  Choose a suitable basis $[P_0],\dots,[P_s]$.
We shall take the \ll standard monomials\rr 
in $h\b_1,\dots,h\b_{n-1}$ with respect to the
reduced Gr\"obner basis for the ideal $(D_1^h,\dots,D_{n-1}^h)$,
using the graded reverse lexicographic monomial order in which
$\b_1,\dots,\b_{n-1}$ are assigned weight one with
$\b_1>\dots>\b_{n-1}$.  As mentioned earlier, the 
differential operators $D_1^h,\dots,D_{n-1}^h$
are the \ll conserved quantities of the quantum Toda lattice\rrr, and
we shall need their explicit form.  They are
obtained from $\Cal R_1,\dots,\Cal R_{n-1}$ simply by replacing $b_i$ by $h\b_i$,
where the relations $\Cal R_i$ are obtained as follows. Let
$$
Z=
\pmatrix
x_1 & q_1     &         &   &        &        &    
\\
-1       & x_2 & q_2     &        &        &        &       
\\
        & -1       & x_3 & q_3    &             &        &       
\\
   &   & \ddots  & \ddots & \ddots      &       &  
\\
        &    &         & -1      & x_{n-2} & q_{n-2}     &       
\\
        &    &    &        & -1           & x_{n-1} &
q_{n-1}\\
        &    &    &   &             & -1           & x_n
\endpmatrix
$$
where $b_i=x_1+\dots+x_i$ $(1\le i\le n-1)$ and $x_1+\dots+x_n=0$.  Then
$
\det(Z+\la I)   
=  \sum_{i=0}^n \Cal R_i \la^i
$
(with $\Cal R_0=0$ and $\Cal R_n=1$).  The basis turns out to be given by
the monomials $(h\b_1)^{i_1}(h\b_2)^{i_2}\dots(h\b_{n-1})^{i_{n-1}}$
with $0\le i_1\le 1$,  $0\le i_2\le 2$ , \dots, $0\le i_{n-1}\le n-1$, ordered
in the above fashion, i.e.\ 
$1< h\b_{n-1}<h\b_{n-2}<\dots< h^m \b_{n-1}^{n-1}\b_{n-2}^{n-2}\dots
\b_{2}^2\b_{1}$.

We compute $\Om^h$
using the Ore algebra package in Maple (\cite{Ma}).  The command
{\it normalf} computes the
product $\b_i P_j$ mod the ideal
$(D_1^h,\dots,D_{n-1}^h)$, and from this one can read off the coefficients
relative to the basis $[P_0],\dots,[P_s]$, and hence the matrix $\Om_i^h$.

Finally we solve the p.d.e.\  system.  Since the equations for $Q_1,\dots,Q_{m-2}$
do not involve $Q_0$, we find these first and then solve 
$Q_0^{-1}dQ_0 = \thzero + [Q_1,\om]$ 
for $Q_0$.  Maple is used again here
in view of the large number of equations.

\no{\it Step 3:}   From $Q_0$ we obtain the new basis $[\hat c_0],\dots,[\hat c_s]$
and hence the quantum evaluation map and quantum product.

We shall illustrate this procedure 
for $n=2,3,4$.  The fact that
$L_+$ is a polynomial function of $q_1,\dots,q_{n-1},h$ allows us to assume
that $Q_i^{[j]} =0$ if $j-i$ is odd (by the argument of Proposition 3.1), and to impose
the basepoint condition $L_+\vert_{q=0} = I$   i.e.\ 
$$
Q_0\vert_{q=0} = I\quad \text{and}\quad Q_1\vert_{q=0} = \dots = Q_{m-2}\vert_{q=0} = 0.
$$

\subhead The case $n=2$
\endsubhead

We have $r=1$, $m=1$, and $s_0=1,s_1=1$.

\no{\it Step 1:}   We have immediately $L_+ = Q_0$ and $Q_0=I$, so there is no differential
equation to solve.  

\no{\it Step 2:}  The ideal is generated by $D^h_1 = h^2\b_1^2$, and the 
basis of $D^h/(D_1^h)$ is $[1],[h\b_1]$.  
We obtain $\Om^h = \frac1h \om$,
i.e.\  $\thi=0$ for all $i$, so no change of basis is needed.

\no{\it Step 3:} The quantum product is determined entirely by 
the relation\footnote{Strictly speaking, this relation should be written
$[[b_1]]\circ[[b_1]]=[[q_1]]$, but we shall omit such brackets in this
section, where no confusion is likely.}
$b_1\circ b_1 = q_1$ of $QH^\ast G/B = QH^\ast \C P^1$.  (Note that, for
any $n$,  the 
formulae $1\circ 1=1$ and $1\circ b_i = b_i\circ 1 = b_i$ follow from
the fact that $1$ is the identity element of $\C[q_1,\dots,q_{n-1}]$.)

\subhead The case $n=3$
\endsubhead

We have $r=2$, $m=3$, and $s_0=1,s_1=2,s_2=2,s_3=1$.  We are dealing
with $(s+1)\times (s+1)$ matrices, where $s+1=s_0+s_1+s_2+s_3=6$; these are
regarded as $(m+1)\times (m+1)$ block matrices, where $m+1=4$.

\no{\it Step 1:} $L_+ = Q_0(I + hQ_1)$ with
$$
\align
Q_0 &= I + Q_0^{[2]} \\
Q_1 &= Q_1^{[3]}.
\endalign
$$
We will have to solve 
$$
dQ_1=\thone 
+ [Q_1,\thzero]  - [Q_1,\om]Q_1
$$
which reduces to
$$
dQ_1^{[3]}=\theta^{(1),[3]}
+ [Q_1,\thzero]^{[3]}  - ([Q_1,\om]Q_1)^{[3]}
=\theta^{(1),[3]}.
$$
Then we will substitute the result in $Q_0^{-1}dQ_0 = \thzero + [Q_1,\om]$
and solve for $Q_0$.

\no{\it Step 2:} 
The ideal is generated by 
$$
D^h_1= h^2\b_1^2 + h^2\b_2^2 - h^2\b_1\b_2-q_1-q_2,\quad
D^h_2= h^3\b_1\b_2^2-h^3\b_1^2 \b_2 + q_1 h\b_2  - q_2 h\b_1
$$
and the basis of $D^h/(D_1^h,D_2^h)$ is given by 
$
1, h\b_2, h\b_1, h^2\b_2^2, h^2\b_2\b_1, h^3\b_2^2\b_1.
$
We have to compute
$\Om^h = \frac1h\om + \thzero + h\thone$ where
$$
\align
\om&= \om^{[-1]} + \om^{[1]} +\om^{[3]} \\
\thzero &= \theta^{(0),[2]} \\
\thone &= \theta^{(1),[3]}.
\endalign
$$
It turns out that $\thone = 0$ and
$$
\om=
\pmatrix
0&0&q_1+q_2&0&0&q_1q_2+q_2^2\\
0&0&0&0&q_1&0\\
1&0&0&0&-q_2&0\\
0&0&-1&0&0&0\\
0&1&1&0&0&0\\
0&0&0&1&1&0
\endpmatrix
dt_1
+
\pmatrix
0&0&0&0&0&q_1q_2+q_2^2\\
1&0&0&q_2&0&0\\
0&0&0&q_2&0&0\\
0&1&0&0&0&-q_2\\
0&0&1&0&0&2q_2\\
0&0&0&0&1&0
\endpmatrix
dt_2
$$
{}
$$
\thzero=
\pmatrix
0&0&0&0&0&0\\
0&0&0&0&0&0\\
0&0&0&0&0&0\\
0&0&0&0&0&0\\
0&0&0&0&0&0\\
0&0&0&0&0&0
\endpmatrix
dt_1
+
\pmatrix
0&0&0&q_2&0&0\\
0&0&0&0&0&0\\
0&0&0&0&0&q_2\\
0&0&0&0&0&0\\
0&0&0&0&0&0\\
0&0&0&0&0&0
\endpmatrix
dt_2.
$$
So $Q_1=0$ and the system reduces to the equation for $Q_0=I + Q_0^{[2]}$:
$$
dQ_0^{[2]} = (Q_0\thzero)^{[2]} = \theta^{(0),[2]}.
$$
By inspection we see that the solution is $Q_0^{[2]}=\theta_2^{(0),[2]}$, so
$$
Q_0=I + \theta_2^{(0),[2]} =
\pmatrix
1&0&0&q_2&0&0\\
0&1&0&0&0&0\\
0&0&1&0&0&q_2\\
0&0&0&1&0&0\\
0&0&0&0&1&0\\
0&0&0&0&0&1
\endpmatrix
$$

\no{\it Step 3:} The quantum product is determined by the quantum evaluation map 
$\hat c_i^\circ = \hat\de[\hat c_i] = [[\hat c_i\vert_{q=0}]]=[[c_i]]$
(and the relations $\Cal R_1,\Cal R_2$).  To calculate this, recall that
$$
c_0=1,\quad
c_1= b_2,\quad
c_2=b_1,\quad
c_3=b_2^2,\quad
c_4=b_2b_1,\quad
c_5=b_2^2b_1.
$$ 
Then
$$
\hat c_0=1,\quad
\hat c_1= b_2,\quad
\hat c_2=b_1,\quad
\hat c_3=b_2^2-q_2,\quad
\hat c_4=b_2b_1,\quad
\hat c_5=b_2^2b_1-q_2 b_1
$$
is obtained by applying the matrix $Q_0^{-1}=I - \theta_2^{(0),[2]}$.
The quantum evaluation map is therefore
$$
b_2\circ b_2 - q_2 = b_2^2,\quad
b_2\circ b_1  = b_2b_1,\quad
b_2\circ b_2\circ b_1 - q_2 b_1 = b_2^2b_1.
$$
These quantum products, together with the relations $\Cal R_1,\Cal R_2$,
determine all other quantum products.

\subhead The case $n=4$
\endsubhead

We have $r=3$, $m=6$, and $s_0=1,s_1=3,s_2=5,s_3=6,s_4=5,s_5=3,s_6=1$. We are
dealing with $24\times 24$ matrices, regarded as $7\times 7$ block matrices.

\no{\it Step 1:} $L_+ = Q_0(I + hQ_1+h^2Q_2+h^3Q_3+h^4Q_4)$ with
$$
\align
Q_0 &= I + Q_0^{[2]} + Q_0^{[4]}+ Q_0^{[6]}\\
Q_1 &= Q_1^{[3]} + Q_1^{[5]}\\
Q_2 &= Q_2^{[4]} + Q_2^{[6]}\\
Q_3 &= Q_3^{[5]} \\
Q_4 &= Q_4^{[6]} .
\endalign
$$
We will have to solve 
$$
\align
dQ_1^{[3]}&=\theta^{(1),[3]}  +  [Q_2^{[4]},\om^{[-1]}]
\\
dQ_1^{[5]}&=\theta^{(1),[5]}  +  [Q_1^{[3]},\theta^{(0),[2]}] 
+ [Q_2^{[4]},\om^{[1]}] + [Q_2^{[6]},\om^{[-1]}]
 - [Q_1^{[3]},\om^{[-1]}]Q_1^{[3]}
\\
dQ_2^{[4]}&=\theta^{(2),[4]}  +  [Q_3^{[5]},\om^{[-1]}]
\\
dQ_2^{[6]}&=\theta^{(2),[6]}  +  [Q_1^{[3]},\theta^{(1),[3]}] 
+ [Q_2^{[4]},\theta^{(0),[2]}]
+ [Q_3^{[5]},\om^{[1]}] 
 - [Q_1^{[3]},\om^{[-1]}]Q_2^{[4]}
\\
dQ_3^{[5]}&=\theta^{(3),[5]}  +  [Q_4^{[6]},\om^{[-1]}]
\\
dQ_4^{[6]}&=\theta^{(4),[6]}
\endalign
$$
The ordering of the matrices here is: 
$Q_1^{[5]} > Q_2^{[6]} > Q_1^{[3]} > Q_2^{[4]} > Q_3^{[5]} > Q_4^{[6]}$.
Beginning with $Q_4^{[6]}$, it is possible to integrate the equations
successively.
Then we will substitute $Q_1$ into $Q_0^{-1}dQ_0 = \thzero + [Q_1,\om]$
and solve for $Q_0$ in a similar manner.

\no{\it Step 2:} 
It turns out that
$\Om^h = \frac1h\om + \thzero + h\thone$, i.e.\ 
$\thtwo = \ththree = \thfour = 0$.  This immediately tells
us that $Q_4^{[6]} = Q_3^{[5]} = Q_2^{[4]} = 0$ and we are left with
the system
$$
\align
dQ_1^{[3]}&=\theta^{(1),[3]}
\\
dQ_1^{[5]}&=\theta^{(1),[5]}  +  [Q_1^{[3]},\theta^{(0),[2]}] 
+ [Q_2^{[6]},\om^{[-1]}]
 - [Q_1^{[3]},\om^{[-1]}]Q_1^{[3]}
\\
dQ_2^{[6]}&=  [Q_1^{[3]},\theta^{(1),[3]}]
\endalign
$$
with $Q_1^{[5]} > Q_2^{[6]} > Q_1^{[3]}$.  This is a system of
$42+1+6=49$ linear equations for the $49$ unknown functions 
in $Q_1,Q_2$ and Maple gives the solution, which turns out to
be: $Q_1= \th^{(1),[3]}$, $Q_2=0$.

A similar calculation can be done for the $1+19+71=91$  unknown functions 
in $Q_0$.   The matrix $Q_0^{-1}$ is given in the Appendix.

\no{\it Step 3:}  The columns of this matrix, interpreted as vectors
with respect to the basis given by the polynomials
$$
\align
&1;\ 
b_3, b_2, b_1;\ 
b_3^2, b_3b_2, b_3b_1, b_2^2, b_2b_1;\ 
b_3^3, b_3^2b_2, b_3^2b_1, b_3b_2^2, b_3b_2b_1, b_2^2b_1;\ 
\\
&\quad \quad 
b_3^3b_2, b_3^3b_1, b_3^2b_2^2, b_3^2b_2b_1, b_3b_2^2b_1;\ 
b_3^3b_2^2, b_3^3b_2b_1, b_3^2b_2^2b_1;\ 
b_3^3b_2^2b_1,
\endalign
$$
produce the quantum evaluation map:

\no
$
(b_3^2 - q_3)^\circ  = b_3^2 
$
\newline
$
(b_3 b_2)^\circ  = b_3b_2
$
\newline
$
(b_3b_1)^\circ  = b_3b_1  
$
\newline
$
(b_2^2 - q_2)^\circ  =b_2^2  
$
\newline
$
(b_2 b_1)^\circ  = b_2b_1 
$

\no
$
(b_3^3 - q_3b_3 - q_3 b_2)^\circ  =b_3^3  
$
\newline
$
(b_3^2 b_2 - q_3b_2)^\circ  =b_3^2b_2  
$
\newline
$
(b_3^2 b_1 - q_3 b_1)^\circ   = b_3^2b_1  
$
\newline
$
(b_3b_2^2 - q_2 b_3)^\circ  = b_3b_2^2  
$
\newline
$
(b_3 b_2 b_1)^\circ  = b_3b_2b_1  
$
\newline
$
(b_2^2 b_1 -q_2 b_1)^\circ  = b_2^2b_1 
$

\no
$
(b_3^3 b_2 - q_3 b_3 b_2 - q_3 b_2^2)^\circ   = b_3^3b_2  
$
\newline
$
(b_3^3 b_1 - q_3 b_3 b_1 - q_3 b_2 b_1)^\circ   =
b_3^3b_1  
$
\newline
$
(b_3^2 b_2^2 - q_2 b_3^2 - q_3 b_2^2)^\circ   =
b_3^2b_2^2  
$
\newline
$
(b_3^2 b_2b_1 - q_3 b_2 b_1)^\circ    =
b_3^2b_2b_1  
$
\newline
$
(b_3 b_2^2 b_1 - q_2 b_3 b_1)^\circ    =
b_3b_2^2b_1 
$

\no
$
(b_3^3 b_2^2
+ 2 q_2q_3 b_3 - 2q_3^2 b_2 - 2 q_2q_3b_1 
 - q_2 b_3^3
+ 2 q_3 b_3^2 b_2 - 3 q_3 b_3 b_2^2)^\circ  =
b_3^3b_2^2  
$
\newline
$
(b_3^3 b_2 b_1
-  q_3 b_3  b_2  b_1 - q_3 b_2^2 b_1)^\circ
 =
b_3^3b_2b_1  
$
\newline
$
(b_3^2 b_2^2 b_1
-  q_2 b_3^2 b_1 - q_3 b_2^2 b_1)^\circ
 =
b_3^2b_2^2b_1 
$

\no
$
(b_3^3 b_2^2b_1
-2 q_1q_2 q_3 - 2 q_2 q_3^2 - 2 q_2^2 q_3
+ 2 q_2q_3 b_3^2 
- 2q_2 q_3 b_3 b_2 
+ 2 q_2q_3 b_3 b_1  + 2 q_2 q_3 b_2^2
-2q_2 q_3 b_2 b_1 - 2 q_3^2 b_2  b_1 
- q_2 b_3^3 b_1
+ 2q_3 b_3^2 b_2  b_1 - 3 q_3 b_3 b_2^2 b_1)^\circ
 =
b_3^3b_2^2b_1 
$

\no These quantum products, together with the relations $\Cal R_1,\Cal R_2,\Cal R_3$,
determine all other quantum products.

\head
\S 4 Quantum Schubert polynomials
\endhead

So far we have used, purely for computational convenience, a monomial basis
$[c_0],\dots,[c_s]$ of $H^\ast G/B$.  The quantum evaluation map $\hat\de$
expresses these basis elements as \ll quantum polynomials\rrr.  Our method
allows us to deduce analogous results for {\it any} basis of $H^\ast G/B$:

\proclaim{Proposition 4.1} 
Let $c'_i = \sum_{j=0}^s C_{ji} c_j$,  for some (constant) matrix $C$
in $\glsc$.
Let $\hat c'_i = \sum_{j=0}^s (Q'_0)^{-1}_{ji} c'_j$
where the matrix function $Q'_0$ is obtained by the method of \S 1, \S 2.
Then $Q'_0= C^{-1} Q_0 C$.
\endproclaim

\demo{Proof}  Let $P'_i = \sum_{j=0}^s C_{ji} P_j$. 
As in \S 1 we obtain the following objects:
$\Om'$ (dropping the superscript $h$ now),
$L'$, $L'_+ = Q'_0 + O(h)$, and $\hat P'_i = \sum_{j=0}^s (L'_+)^{-1}_{ji} P'_j$.
Let us compare $\Om'$ with $\Om$  ($=\Om^h$ in our earlier notation).
We have 
$\b_j P'_i = \sum_k C_{ki} \b_j P_k = \sum_{k,l} C_{ki} (\Om_j)_{lk} P_l$
and
$\b_j P'_i = \sum_{k} (\Om'_j)_{ki} P'_k  = \sum_{k,l} (\Om'_j)_{ki} C_{lk} P_l$,
hence $\Om_j C = C \Om'_j$ and $\Om' = C^{-1} \Om C$.  It follows that
$L' = XLC$ for some constant invertible matrix $X$, hence
$L'_- L'_+ = (X L_- X^{-1}) (X L_+ C)$ and $L'_+ = XL_+ C$.  The
basepoint condition $L'_+ \vert_{q=0} = I$ gives $X = C^{-1}$.  So
$L'_+ = C^{-1}L_+C$ and $Q_0' = C^{-1} Q_0 C$.
\qed\enddemo

\proclaim{Corollary 4.2}  We have
$\hat c'_i = \sum_{k=0}^s R_{ki} c_k$ where $R = C (Q_0' )^{-1} = Q_0^{-1} C$.
\endproclaim

\demo{Proof} The formula $R = C (Q_0' )^{-1}$
is immediate from $\hat c'_i = \sum_{j=0}^s
(Q'_0)^{-1}_{ji} c'_j$ and $c'_i = \sum_{j=0}^s C_{ji} c_j$.
Hence  $R  = Q_0^{-1} C$ by Proposition 4.1.
\qed\enddemo

As an example, take $c'_0,\dots,c'_s$ to be the  Schubert polynomials.  
Then $\hat c'_0,\dots,\hat c'_s$ are {\it quantum} Schubert polynomials
in the sense that they are polynomials whose \ll quantum
evaluations\rr produce exactly the Schubert classes:
$$
(\hat c'_i)^\circ = 
\de [\hat c'_i ] = [[ \hat c'_i\vert_{q=0} ]] = [[ c'_i ]].
$$
Thus, from $C$ (which is well known) and $Q_0$ (which we have calculated), 
we obtain $R$, and hence 
{\it explicit formulae for these  quantum Schubert polynomials.}
Observe that $\hat c'_i\vert_{q=0} = c'_i$, as $Q'_0\vert_{q=0}=I$,
so our quantum Schubert polynomials are indeed  \ll $q$-deformations\rr of the
Schubert polynomials.

In \S 14 of \cite{Fo-Ge-Po}, tables of quantum Schubert polynomials (hence
$R$ and $C$) are given for $G=\glnc$ and $n=2,3,4$. It is easy to verify
that these coincide with ours (and it is natural to conjecture that
this holds for general $n$).
To convert from \cite{Fo-Ge-Po} to our
notation, $x_1,\dots,x_{n-1}$ should be replaced by 
$b_{n-1}, b_{n-2}-b_{n-1}, \dots, b_{1}-b_{2}$ and 
$q_1,\dots,q_{n-1}$ by $q_{n-1},\dots,q_{1}$. 
For example, when $n=2$, the Schubert polynomials from \cite{Fo-Ge-Po}
are
$$
1, b_2, b_1, -b_2^2+b_2b_1, b_2^2, b_2^2b_1,
$$
so (with respect to the usual 
monomial basis given by
$1,b_2,b_1,b_2^2,b_2b_1,b_2^2b_1$)
we have
$$
C =
\pmatrix
1 & 0 & 0 & 0 & 0 & 0 \\
0 & 1 & 0 & 0 & 0 & 0  \\
0 & 0 & 1 & 0 & 0 & 0  \\
0 & 0 & 0 & -1 & 1 & 0  \\
0 & 0 & 0 & 1 & 0 & 0  \\
0 & 0 & 0 & 0 & 0  & 1 
\endpmatrix.
$$
Using our matrix $Q_0$ from \S 3, the matrix $R= Q_0^{-1} C$ is
$$
\pmatrix
1&0&0&-q_2&0&0\\
0&1&0&0&0&0\\
0&0&1&0&0&-q_2\\
0&0&0&1&0&0\\
0&0&0&0&1&0\\
0&0&0&0&0&1
\endpmatrix
\pmatrix
1 & 0 & 0 & 0 & 0 & 0  \\
0 & 1 & 0 & 0 & 0 & 0  \\
0 & 0 & 1 & 0 & 0 & 0  \\
0 & 0 & 0 & -1 & 1 & 0  \\
0 & 0 & 0 & 1 & 0 & 0 \\
0 & 0 & 0 & 0 & 0 & 1 
\endpmatrix
=
\pmatrix
1 & 0 & 0 & q_2 & -q_2 & 0 \\
0 & 1 & 0 & 0 & 0 & 0  \\
0 & 0 & 1 & 0 & 0 & -q_2  \\
0 & 0 & 0 & -1 & 1 & 0  \\
0 & 0 & 0 & 1 & 0 & 0 \\
0 & 0 & 0 & 0 & 0  & 1 
\endpmatrix.
$$
Hence our quantum Schubert polynomials are 
$$
1, b_2, b_1, -b_2^2+b_2b_1+q_2, b_2^2-q_2, b_2^2b_1-q_2b_1,
$$
in agreement with those of \cite{Fo-Ge-Po}.  
The case $n=4$ may be verified
in the same way, by reading off $C$ from \cite{Fo-Ge-Po} and
computing $R=Q_0^{-1} C$, where $Q_0^{-1}$ is given in the
Appendix.

\newpage
\head
Appendix: $Q_0^{-1}$ for $G=GL_4\C$
\endhead

{
%\magnification=1000
\eightpoint
$Q_0^{-1}=
\left(
\matrix
1 &  0 &  0 &  0 &  -q_3 &  0 &  0 &  -q_2 &  0 &  
0 &  0 &  0 &  0 &  0 &  0 
\\
0 &  1 &  0 &  0 &  0 &  0 &  0 &  0 &  0 & 
-q_3 &  0 &  0 &  -q_2 &  0 &  0 
\\
0 &  0 &  1 &  0 &  0 &  0 &  0 &  0 &  0 &  -q_3 &  -q_3 &  0 &  0 &  0 &  0 
\\
0 &  0 &  0 &  1 &  0 &  0 &  0 &  0 &  0 &  0 &  0 &  -q_3 &  0 &  0 &  -q_2 
\\
0 &  0 &  0 &  0 &  1 &  0 &  0 &  0 &  0 &  0 &  0 &  0 &  0 & 
 0 &  0 
\\
0 &  0 &  0 &  0 &  0 &  1 &  0 &  0 &  0 &  0 & 
0 &  0 &  0 &  0 &  0
\\ 
0 &  0 &  0 &  0 &  0 &  0 &  1 &  0 &  0 &  0 &  0 &  0 &  0 &  0 &  0 
\\
0 &  0 &  0 &  0 &  0 &  0 &  0 &  1 &  0 &  
0 &  0 &  0 &  0 &  0 &  0
\\ 
0 &  0 &  
0 &  0 &  0 &  0 &  0 &  0 &  1 &  0 &  0 &  0 &  0 &  0 &  0 
\\
0 &  0 &  0 &  0 &  0 &  0 & 
0 &  0 &  0 &  1 &  0 &  0 &  0 &  0 &  0 
\\
0 &  0 &  0 &  0 &  0 &  0 &  0 &  0 &  0 &  0 &  1 &  0 &  0 &  0 &  0 
\\
0 &  0 &  0 &  0 &  0 &  0 &  0 &  0 &  0 &  0 &  0 &  1 &  0 & 
 0 &  0 
\\
0 &  0 &  0 &  0 &  0 &  0 &  0 &  0 & 
0 &  0 &  0 &  0 &  1 &  0 &  0 
\\
0 &  0 &  0 &  0 &  0 &  0 &  0 &  0 &  0 &  0 &  0 &  0 &  0 &  1 &  0 
\\
0 &  0 &  0 &  0 &  0 &  0 &  0 &  0 &  0 &  0 &  0 &  0 &  0 &  0 & 
 1 
\\
0 &  0 &  0 &  0 &  0 &  0 &  0 &  0 & 
0 &  0 &  0 &  0 &  0 &  0 &  0 
\\
0 &  0 &  0 &  0 &  0 &  0 &  0 &  0 &  0 &  0 &  0 &  0 &  0 &  0 &  0 
\\
0 &  0 &  0 &  0 &  0 &  0 &  0 &  0 &  0 &  0 &  0 &  0 &  0 &  0 &  0 
\\
0 &  0 &  0 &  0 &  0 &  0 &  0 &  0 &  0 &  0 &  0 &  0 &  0 & 
0 &  0
\\
0 &  0 &  0 &  0 &  0 &  0 &
  0 &  0 &  0 &  0 &  0 &  0 &  0 &  0 &  0
\\
0 & 
0 &  0 &  0 &  0 &  0 &  0 &  0 &  0 &  0 &  0 &  0 &  0 &  0 &  0 
\\
0 &  0 &  0 &  0 &  0 &  0 &  0 &  0 &  0 &  0 &  0 &  0 &  0 &  0 &  0 
\\
0 &  0 &  0 &  0 &  0 &  0 &  0 &  0 &  0 &  0 & 
 0 &  0 &  0 &  0 &  0 
\\
0 &  0 &  0 &  0 &  0 &  0 &  0 &  0 & 
0 &  0 &  0 &  0 &  0 &  0 &  0 
\endmatrix
\right.
$
}

{
%\magnification=1000
\eightpoint
$
\left.
\matrix
  0 &  0 &  0 &  0 &  0 &  0 &  0 &  0 &  -2 q_1 q_2 q_3-2 q_2 q_3^2-2 q_2^2 q_3
\\
  0 &  0 &  0 &  0 &  0 &  2 q_2 q_3 &  0 &  0 &  0
\\
  0 &  0 &  0 &  0 &  0 &  -2 q_3^2 &  0 &  0 &  0
\\
  0 &  
0 &  0 &  0 &  0 &  -2 q_2 q_3 &  0 &  0 & 
0
\\
  0 &  0 &  -q_2 &  0 &  0 &  0 &  0 &  0 &  2 q_2 q_3
\\
  -q_3 &  0 &  0 &  0 &  0 &  0 &  0 &  0 &  -2 q_2 q_3
\\
  0 &  -q_3 &  0 &  0 & 
-q_2 &  0 &  0 &  0 &  2 q_2 q_3
\\
  -q_3 &  0 &  -q_3 &  0 &  0 &  0 &  0 &  0 &  2 q_2 q_3
\\ 
  0 & 
 -q_3 &  0 &  -q_3 &  0 &  0 &  0 &  0 &  -2 q_2 q_3-2 q_3^2
\\
  0 &  0 &  0 &  0 &  0 &  -q_2 &  0 &  0 &  0
\\
  0 &  0 &  0 &  0 &  0 & 
2 q_3 &  0 &  0 &  0
\\
  0 &  0 &  0 &  0 &  0 &  0 &  0 &  -q_2 &  0
\\
  0 &  0 &  0 &  0 &  0 &  -3 q_3 &  0 &  0 &  0
\\
  0 &  0 &  0 &  0 &  0 &  0 & 
-q_3 &  0 &  0
\\
  0 &  0 &  0 &  0 &  0 &  0 &  -q_3 &  -q_3 &  0
\\
  1 &  0 &  0 &  0 &  0 &  0 &  0 &  0 &  0
\\
 0 &  1 & 
 0 &  0 &  0 &  0 &  0 &  0 & 
-q_2
\\
  0 & 
 0 &  1 &  0 &  0 &  0 &  0 &  0 &  0
\\
  0 &  0 &  0 &  1 &  0 &  0 &  0 &  0 &  2 q_3
\\
   0 &  0 &  0 &  0 &  1 &  0 &  0 &  0 &  -3 q_3
\\
  0 &  0 &  0
 &  0 &  0 &  1 &  0 &  0 &  0
\\
  0 &  0 & 
0 &  0 &  0 &  0 &  1 &  0 &  0
\\
  0 &  0 &  0 &  0 &  0 &  0 &  0 &  1 &  0
\\
  0 &  0 &  0 &  0 &  0 &  0 &  0 &  0 &  1 
\endmatrix
\right)
$
}

\bigskip
\it
\no  
Department of Mathematics,
Graduate School of Science,
Tokyo Metropolitan University,
Minami-Ohsawa 1-1, Hachioji-shi,
Tokyo 192-0397, Japan
\medskip
Current address of first author:
\smallskip
\no School of Mathematics and Computer Science,
National University of Mongolia,
Ulannbaatar, Post Office 46-A, Post Box 586, Mongolia


\newpage

\parskip 0pt

\Refs   

\eightpoint

\widestnumber\key{XXXXXX}

\ref
\key Bu-Gu
\by F.E. Burstall and M.A. Guest
\paper Harmonic two-spheres in compact symmetric spaces, revisited
\jour Math. Ann. 
\vol 309
\yr 1997
\pages 541--572 
\endref

\ref
\key Ci
\by I. Ciocan-Fontanine
\paper Quantum cohomology of flag varieties
\jour Int. Math. Res. Notices
\vol 6
\yr 1995
\pages 263--277
\endref

\ref
\key Fo-Ge-Po
\by S. Fomin, S. Gelfand, and A. Postnikov
\paper Quantum Schubert polynomials
\jour J. Amer. Math. Soc. 
\vol 10 
\yr 1997
\pages 565--596
\endref

\ref
\key  Gi
\by A.B. Givental
\paper Equivariant Gromov-Witten invariants
\jour Internat. Math. Res. Notices
\yr 1996
\vol 13
\pages 1--63
\endref

\ref 
\key  Gi-Ki
\by  A. Givental and B. Kim
\paper Quantum cohomology of flag manifolds and Toda lattices
\jour Commun. Math. Phys.
\yr 1995
\vol 168
\pages 609--641
\endref

\ref
\key Go-Wa
\by R. Goodman and N. Wallach
\paper Classical and quantum mechanical systems of Toda-lattice type. III. 
Joint eigenfunctions of the quantized systems
\jour Commun. Math. Phys.
\yr 1986
\vol 105
\pages 473--509
\endref

\ref
\key  Gu
\by  M.A. Guest
\paper Quantum cohomology via D-modules
\jour Topology (math.DG/0206212)
\yr 
\vol 
\pages 
\paperinfo to appear
\endref

\ref
\key  Ki
\by B. Kim
\paper Quantum cohomology of flag manifolds $G/B$ and
quantum Toda lattices
\jour Ann. of Math.
\yr 1999
\vol 149
\pages 129--148
\endref
%alg-geom/9607001

\ref
\key Ki-Ma
\by A.N. Kirillov and T. Maeno
\paper Quantum double Schubert polynomials, quantum Schubert 
polynomials and Vafa-Intriligator formula
\jour Discrete Math.
\vol 217 
\yr 2000
\pages 191--223
\endref

\ref
\key Ko
\by B. Kostant
\paper The solution to a generalized Toda lattice
and representation theory
\jour Adv. Math.
\vol 34
\yr 1979
\pages 195--338
\endref

\ref
\key Ma
\by 
\paper Maple 7
\jour Waterloo Maple
\vol 
\yr 1981--2001
\pages 
\endref

\ref
\key  Pr-Se
\by A.N. Pressley and G.B. Segal
\book  Loop Groups
\publ Oxford Univ. Press 
\yr 1986
\endref

\endRefs
\enddocument